\documentclass[oneside]{amsart}
\usepackage{enumerate,amssymb}
\usepackage{graphicx}
\newtheorem{thm}{Theorem}[section]
\newtheorem*{thmn}{Theorem}
\newtheorem{coro}[thm]{Corollary}
\newtheorem{prop}[thm]{Proposition}

\newtheorem*{claim}{Claim}

\newtheorem{lem}[thm]{Lemma}
\theoremstyle{definition}
\newtheorem{defn}[thm]{Definition}

\newtheorem{question}[thm]{Question}
\newcommand{\Rset}{\mathbb{R}}
\newcommand{\Cset}{2^\omega}
\newcommand{\cset}{2^{<\omega}}
\newcommand{\Nset}{\mathbb{N}}
\newcommand{\Zset}{\mathbb{Z}}
\newcommand{\Tset}{\mathbb{T}}
\newcommand{\Pset}{\omega^\omega}
\newcommand{\abs}[1]{\lvert#1\rvert}
\newcommand{\norm}[1]{\lVert#1\rVert}
\newcommand{\seq}[1]{\langle#1\rangle}
\newcommand{\hm}{\mathcal H}
\newcommand{\eps}{\varepsilon}
\newcommand{\del}{\delta}
\newcommand{\subs}{\subseteq}
\newcommand{\wh}{\widehat}
\newcommand{\Implies}{$\Rightarrow$}
\renewcommand{\leq}{\leqslant}
\renewcommand{\geq}{\geqslant}
\DeclareMathOperator{\ord}{ord}
\newcommand{\leb}{\mathcal L}
\DeclareMathOperator{\hdim}{\dim_{\mathsf{H}}}
\DeclareMathOperator{\pdim}{\dim_{\mathsf{P}}}
\DeclareMathOperator{\diam}{diam}
\DeclareMathOperator{\dist}{dist}
\newcommand{\lowerd}{\underline{d}}
\newenvironment{enum}{\begin{enumerate}[\rm(i)]}{\end{enumerate}}
\newenvironment{itemyze}%
  {\begin{list}{\textbullet}{\labelwidth1ex\setlength{\leftmargin}{1.5em}}}%
  {\end{list}}
\newcommand{\embed}{\hookrightarrow}
\newcommand{\concat}{^{\mkern-2mu\frown}\mkern-3mu}
\newcommand{\si}{$\sigma$\nobreakdash-}
\begin{document}
\title
[Mapping by Lipschitz and nearly Lipschitz maps]
{Mapping Borel sets onto balls and self-similar sets by Lipschitz and nearly
Lipschitz maps}

\author{Ond\v rej Zindulka}
\address
{Department of Mathematics\\
Faculty of Civil Engineering\\
Czech Technical University\\
Th\'akurova 7\\
160 00 Prague 6\\
Czech Republic}
\email{ondrej.zindulka@cvut.cz}
\urladdr{http://mat.fsv.cvut.cz/zindulka}
\thanks{The author was supported by the grant GA\v CR 15--082185
of the Grant Agency of the Czech Republic.}

\begin{abstract}
If $X$ is an analytic metric space satisfying a very mild doubling condition,
then for any finite Borel measure $\mu$ on $X$ there is a set $N\subs X$
such that $\mu(N)>0$, an ultrametric space $Z$ and a Lipschitz bijection
$\phi:N\to Z$ whose inverse is nearly Lipschitz,
i.e., $\beta$-H\"older for all $\beta<1$.

As an application it is shown that a Borel set in a Euclidean space maps onto $[0,1]^n$
by a nearly Lipschitz map if and only if it cannot be covered by
countably many sets of Hausdorff dimension strictly below $n$.

The argument extends to analytic metric spaces satisfying the mild condition.
Further generalization replaces cubes with self-similar sets,
nearly Lipschitz maps with nearly H\"older maps and integer dimension with arbitrary
finite dimension.
\end{abstract}

\maketitle

\section{Introduction}

It is easy to prove that every compact set $K\subs\Rset$ of real numbers with
positive Lebesque measure can be mapped onto the interval $[0,1]$ by a Lipschitz map.
In~\cite{MR1147388}, Mikl\'os Laczkovich asked if this remains true in higher dimensions.
In more detail, if it is true that for every natural number $n>1$ and
every compact set $K\subs\Rset^n$ with positive $n$-dimensional Lebesgue measure
there is a Lipschitz mapping $f:K\to[0,1]^n$ onto the $n$-dimensional cube.
This question turned to be very difficult. So far the (affirmative) answer was found only
for $n=2$ (by David Preiss, published years later in~\cite{MR2185733}).

Vitushkin, Ivanov and Melnikov~\cite{MR0154965} (see also~\cite{MR1313694})
constructed an example that shows that no generalization beyond the Laczkovich's
question is possible: a compact set $K\subs\Rset^2$ with positive linear measure
(i.e., the $1$-dimensional Hausdorff measure)
that cannot be mapped onto a segment by a Lipschitz map.

Let $n$ be a positive natural number and let us denote
the $n$-dimensional Hausdorff measure by $\hm^n$ and the Hausdorff dimension by $\hdim$.
Using recent results on \emph{monotone metric spaces}~\cite{MR2957686}
and ultrametric spaces~\cite{MR3032324} (see below)
Keleti, M\'ath\'e and Zindulka~\cite{MR3159074} proved that if the assumption
$\hm^n(K)>0$ is strengthened to $\hdim K>n$,
then $K$ can be mapped by a Lipschitz map onto $[0,1]^n$ for any analytic metric space $K$.

From what have been said, it is clear that we do not have a complete understanding of
\begin{itemyze}
\item what condition akin to $\hm^n(K)>0$ is equivalent to the existence of
a Lipschitz mapping of $K$ onto $[0,1]^n$,
\item what remains true about mapping of $K$ onto $[0,1]^n$ if we
merely suppose that $\hm^n(K)>0$.
\end{itemyze}
The present paper deals with the latter question and provides a partial answer.

We build upon ideas presented in \cite{MR3032324}
and \cite{MR3159074}. In~\cite{MR3032324} Mendel and Naor prove that a compact metric
space contains a large (with respect to Hausdorff dimension) subset with a rather
simple metric structure -- it is, up to bi-Lipschitz homeomorphism, ultrametric.
We also construct a large subset with a simple structure,
the ``large'' and ``simple'' notions, however, are a bit different.

In more detail, we prove in Section~\ref{sec:ultra} a theorem
about spaces satisfying a mild doubling condition henceforth termed
\emph{non-exploding} spaces.
Roughly speaking, the condition requires that the number of balls of radius $r$ needed
to cover $X$ does not increase too fast with decreasing $r$, cf.~Definition~\ref{Nee}.
The following is a simplified version of Theorem~\ref{main1}.
\begin{thmn}
Let $X$ be a non-exploding analytic metric space
and let $\mu$ be a finite Borel measure on $X$. Then there is a set $N\subs X$
such that $\mu(N)>0$, an ultrametric space $Z$ and a Lipschitz bijection
$\phi:N\to Z$ whose inverse is $\beta$-H\"older for all $\beta<1$.
\end{thmn}

Then, based upon this result, we show
that if a non-exploding analytic metric space
has positive $n$-dimensional Hausdorff measure, then
there is a mapping of $X$ onto $[0,1]^n$ that is as close to Lipschitz as it gets.
(We know from the Vitushkin's example that it may fail to be Lipschitz.)

\begin{thmn}
Let $X$ be an analytic non-exploding space and $n\in\Nset$. If $\hm^n(X)>0$,
then there is a mapping $f:X\to[0,1]^n$ onto $[0,1]^n$ that is
$\beta$-H\"older for every $\beta<1$.
\end{thmn}

This is proved in Section~\ref{sec:maps}. The full strength of the mapping
theorem is stated in Theorem~\ref{mainlip} and its corollaries \ref{corocube}
and ~\ref{coroX}.
%

Once we have this mapping result, we may ask if the sufficient condition
on the Hausdorff measure is also necessary. And it turns out that it is not,
nevertheless there is a simple, natural condition on $X$ involving
Hausdorff dimension that is necessary and also sufficient.
We discuss this in Section~\ref{sec:indec}.

Further generalization replaces cubes with self-similar sets,
nearly Lipschitz maps with nearly H\"older maps and integer dimension with arbitrary
finite dimension.

The last Section~\ref{sec:rem} contains remarks and presents some questions and
problems.
\smallskip

The results of this paper already found an application. Namely,
Balka, Elekes and M\'ath\'e use them in~\cite{Balka2016221} to prove the following theorem
on a prevalent behavior of continuous functions.
(The dimensions involved are the Hausdorff and packing dimensions of sets and measures.)

\begin{thmn}[{\cite{Balka2016221}}]
Let $K$ be a non-exploding compact metric space
and $\mu$ a continuous, finite Borel measure on $K$. Let $n$ be a positive integer.
Then for almost every continuous function on $K$
(in the sense of Christensen's~\cite{MR0326293} Haar measure zero)
with values in $\Rset^n$ there is an open set $U_f\subs\Rset^n$
such that $\mu(f^{-1}(U_f))=\mu(K)$ and for all $y\in U_f$
$$
  \hdim f^{-1}(y)\geq\hdim\mu\text{\quad and \quad}\pdim f^{-1}(y)\geq\pdim\mu.
$$
\end{thmn}
This theorem generalizes a number of previous results and as far as I know, it is
the first theorem of its kind proved in such a general context.
\medskip

All spaces we work with are separable metric spaces. Recall that a metric space
is \emph{analytic} if it is a continuous image of a complete metric space
(or, equivalently, of the irrational numbers, or equivalently, a Suslin set in
a complete metric space). A continuous image of an analytic space is analytic.
Every analytic space is separable.

Some of the common notation includes
$B(x,r)$ for the closed ball centered at $x$, with radius $r$;
$\diam E$ for the diameter of a set $E$ in a metric space;
$\dist(A,B)$ the (lower) distance of two sets $A,B$;
$n,m$ are generic symbols for positive integers.
$\Rset$ and $\Zset$ have the usual meaning;
$\omega$ stands for the set of natural numbers including zero.
$\Cset$ denotes the set of binary sequences; it is also a
compact topological space homeomorphic to the standard Cantor ternary set
and a topological group.
Likewise, $\Pset$ denotes the set of all sequences of natural numbers,
i.e., the maps $f:\omega\to\omega$.

\section{Large nearly ultrametric sets}\label{sec:ultra}

In~\cite{MR3032324}, Mendel and Naor proved that every analytic metric space $X$
contains large ultrametric-like subspaces. In more detail, for every $\eps>0$
there is a subset $Y\subs X$ with the following two properties:
\begin{itemyze}
  \item[(1)] Hausdorff dimension of $Y$ is large: $\hdim Y\geq\hdim X-\eps$,
  \item[(2)] the metric structure of $Y$ is ultrametric-like:
    there is a bijection $f:Y\to U$ onto an ultrametric space $U$ such that both
    $f$ and its inverse are Lipschitz.
\end{itemyze}
Though it is not clear at first glance how this powerful result is related to
mapping metric spaces onto cubes, it is one of two crucial ingredients of the
theorems of~\cite{MR3159074} mentioned in the previous paragraph and its spirit,
as we shall see, is also important for the present paper.

In this section we attempt to prove a theorem similar to that of Mendel and Naor:
our goal is to find, just like in the theorem, within an analytic metric space $X$
a large ultrametric-like set $Y$.
We want a bit more than (1): given a finite Borel measure on $X$,  our set will have
to have positive measure. In order to achieve that, we have to sacrifice some of (2):
our set will be still ultrametric-like, but not quite as much as in (2).
We term the notion \emph{nearly ultrametric};
it is introduced in Definition~\ref{def:nU} below.

We do not succeed completely: our proof only works within a framework of analytic spaces
that are subject to a growth condition similar to the doubling condition, but much weaker,
the so called \emph{non-exploding spaces}.
The result is stated in Theorem~\ref{main1}.

\subsection*{Non-exploding spaces}

We first discuss the doubling-like condition.
Recall that a metric space $X$ is
\emph{doubling} if there is a number $Q$ such that every ball in the space $X$
can be covered by at most $Q$ many balls of halved radii.
This notion and equivalent or similar notions
have been defined, investigated and used throughout the literature.

We may generalize the notion as follows. Suppose that the number $Q$ is not fixed
but may increase as the radius of the ball in question decreases; but it may
not increase very fast.
In more detail:

\begin{defn}\label{Nee}
Let $X$ be a metric space. If there is a function $Q:(0,\infty)\to\Rset$
such that every closed ball in $X$ of radius $r>0$ is covered by at most
$Q(r)$ many closed balls of radius $r/2$ and such that
\begin{equation}\label{def:NEXP}
  \lim_{r\to0}\frac{\log Q(r)}{\log r}=0 \tag{\textsf{NE}},
\end{equation}
we call the metric space $X$ \emph{non-exploding}%
\footnote{The term was coined by Tam\'as Keleti.}.
\end{defn}
Needless to say that every doubling metric space
and in particular every subset of a Euclidean space is non-exploding.

\subsection*{Nearly Lipschitz maps}

We will be frequently making use of the notion of a nearly Lipschitz map
that was introduced in~\cite{MR2957686}. We present two of the several equivalent
definitions.
Recall that, given $\beta>0$, a mapping $f:(X,d_X)\to(Y,d_Y)$ is \emph{$\beta$-H\"older}
if there is an $\eps>0$ and a constant $C$ such that if $d_X(x,y)\leq\eps$, then
$d_Y(f(x),f(y))\leq C\,d_X(x,y)^\beta$. (Note that we deviate slightly
from the usual definition by introducing $\eps$, but if $Y$ is bounded,
the definitions are equivalent and, moreover, since we are interested in low scale
behavior, we may always suppose that $Y$ is bounded.)
\begin{defn}[\cite{MR2957686}]
A mapping $f:(X,d_X)\to(Y,d_Y)$ between metric spaces is termed \emph{nearly Lipschitz}
if $f$ is $\beta$-H\"older for all $\beta<1$.
\end{defn}
\begin{prop}[\cite{MR2957686}]\label{prop:nL}
A mapping $f:(X,d_X)\to(Y,d_Y)$ is nearly Lipschitz if and only if
there is a function $h:(0,\infty)\to(0,\infty)$ such that
$\varliminf_{r\to0}\frac{\log h(r)}{\log r}\geq1$ and
\begin{equation}\label{eq:nL2}
  d_Y\bigl(f(x),f(y)\bigr)\leq h\bigl(d_X(x,y)\bigr),
  \quad x,y\in X.
\end{equation}
\end{prop}
\begin{proof}
Suppose that $f:X\to Y$ is nearly Lipschitz. Then there is a decreasing sequence
$\del_n\to 0$ and a sequence of constants $C_n$ such that
$$
  d_X(x,y)\leq\del_n \implies d_Y(f(x),f(y))\leq C_n d_X(x,y)^{1-1/n}.
$$
We may also suppose that $\del_n\leq C_n^{-n}$ so that $\del_n^{-1/n}$
is an upper estimate of $C_n$.
Define the function $h$ by
$$
  h(r)=r^{1-2/n}\text{ if } r\in[\del_{n+1},\del_n).
$$
Since $\log h(r)/\log r=1-2/n$ on the entire interval $[\del_{n+1},\del_n)$,
we have
$
  \varliminf_{r\to0}\frac{\log h(r)}{\log r}=\varliminf_{n\to\infty}1-\frac2n=1.
$

If $\del_{n+1}\leq d_X(x,y)<\del_n$, then
$$
  d_Y(f(x),f(y))\leq\del_n^{-1/n}d_X(x,y)^{1-1/n}\leq d_X(x,y)^{1-2/n}=h(d_X(x,y))
$$
which proves \eqref{eq:nL2}.

The reverse implication is straightforward.
\end{proof}
We will say that the metric spaces $X,Y$ are \emph{nearly Lipschitz equivalent}
if there is a bijective mapping $f:X\to Y$ such that both $f$ and its inverse are
nearly Lipschitz.

\subsection*{Nearly ultrametric spaces}

Recall that a metric
space $(X,d)$ is \emph{ultrametric} if any triple $x,y,z\in X$ of points
satisfies $d(x,z)\leq\max\{d(x,y),d(y,z)\}$.

\begin{defn}\label{def:nU}
A metric space $X$ is \emph{nearly ultrametric} if it is nearly Lipschitz equivalent
to an ultrametric space.
\end{defn}
\begin{prop}\label{nUlt}
Let $X$ be a metric space. The following are equivalent.
\begin{enum}
\item $X$ is nearly ultrametric,
\item there is a nearly Lipschitz bijection $f:X\to Y$ onto an ultrametric space
with Lipschitz inverse,
\item there is a Lipschitz bijection $f:X\to Y$ onto an ultrametric space
with nearly Lipschitz inverse.
\end{enum}
\end{prop}
\begin{proof}
Of course it is enough to prove (i)\Implies(ii) and (i)\Implies(iii).
Using Proposition~\ref{prop:nL} there are functions $g,h$ such that
$\varliminf_{r\to0}\frac{\log h(r)}{\log r}\geq1$ and
$\varliminf_{r\to0}\frac{\log g(r)}{\log r}\geq1$ and such that
$$
  d_X(x,y)\leq g(d_Y(f(x),f(y)))\leq g\circ h(d_X(x,y)).
$$
(The latter inequality holds if $g$ is non-decreasing, but we may suppose that.)
Define a new metric on $Y$ by $d_Y'=g\circ d_Y$. Since $d_Y$ is an ultrametric, so
is $d_Y'$ and it is easy to verify that the mapping $f:(X,d_X)\to(Y,d_Y')$
is nearly Lipschitz and has Lipschitz inverse. This proves (i)\Implies(ii) and
(i)\Implies(iii) is proved likewise.
\end{proof}

We have enough to state the first summit of the paper.
\begin{thm}\label{main1}
Let $X$ be a non-exploding analytic metric space
and let $\mu$ be a finite Borel measure on $X$. Then for every $\eps>0$ there is
a compact nearly ultrametric set $C\subs X$ such that $\mu(X\setminus C)<\eps$.
\end{thm}

The rest of this section is devoted to the proof of this theorem.
It has two stages: we first prove a very particular case and
then reduce the theorem to this particular case.

\subsection*{The infinite dimensional torus}

We will prove the theorem first for the infinite dimensional torus, i.e.,
the compact group $(\Rset/\Zset)^\omega$ equipped with a special metric,
and its Haar measure.

We will identify $\Rset/\Zset$ with the interval $[0,1)$. The group operation
on $[0,1)$ is of course addition modulo $1$. For $x,y\in[0,1)$ the distance from
$x$ to $y$  is given
by $d(x,y)=\norm{x-y}$, where $\norm z=\min\{\abs z,\abs{1-z}\}$, i.e.,
$[0,1)$ is a circle and $\norm{x-y}$ is the length of the shorter of the two
arcs between $x$ and $y$.
We also consider the Lebesgue measure on $[0,1)$; let us denote it by $\leb$.

Let $\Tset=[0,1)^\omega$ and denote the probability
Haar measure on $\Tset$  by $\lambda_\Tset$
(which is clearly obtained as a product measure from $\leb$).

Now let $G\in\Pset$ be a sequence that is not eventually zero.
Such a function defines a partition of $\omega=I_0\cup I_1\cup I_2\cup\dots$
into disjoint consecutive intervals such that the length of each $I_n$ is $G(n)$.
(Some of the intervals may be empty.)
For each $n$ and $k\in I_n$ let $r_k=2^{-n}$.
Clearly $r_k\searrow 0$.
Now define the metric on $\Tset$ as follows:
$$
  d_G(x,y)=\sup_{k\in\omega}r_k\norm{x_k-y_k}.
$$
It is easy to check that $d_G$ is an invariant metric on $\Tset$.
It is clear that if $G$ assumes large values, then the diameters of the coordinate
circles decrease slowly, and vice versa.
We impose the following growth condition upon $G$. From now on we define $\log 0=0$.

\begin{defn}
We call a sequence $G\in\Pset$ \emph{slow} if
$G(n)>0$ for infinitely many $n$ and $\log G(n)/n\to 0$.
\end{defn}

\begin{lem}\label{lem:haar}
If $G$ is slow, then for every $\eps>0$ there is a compact nearly ultrametric set
$C\subs(\Tset,d_G)$ such that $\lambda_\Tset(C)>1-\eps$.
\end{lem}
\begin{proof}
We first construct the set and then prove that it has the required properties.
Begin with defining two sequences of positive real numbers $\seq{a_m},\seq{b_k}$:
Let
$$
  a_m=\frac{1}{2(m+1)^2},\qquad b_k=\frac{\eps}{4(n+1)^2G(n)},
$$
where $n$ is the unique number such that $k\in I_n$. Clearly
\begin{equation}\label{posl}
  \sum_{m\in\omega}a_m<1,\qquad \sum_{k\in\omega}b_k<\eps/2.
\end{equation}
We now set up, for each $k$, a Cantor-like ternary set $C_k\subs[0,1)$.
Let $\cset=\bigcup_{n\in\omega}2^n$ be the binary tree consisting of all
$\{0,1\}$-valued finite sequences.
Build recursively a family of closed intervals $\{J_s:s\in\cset\}$
as follows. Let $J_\emptyset=[0,1-b_k]$ and if $s\in\cset$ and $J_s$ is constructed,
define $J_{s\concat0}$ and $J_{s\concat1}$ subject to the following conditions:
\begin{enumerate}[\rm(a)]
\item $J_{s\concat0}$ and $J_{s\concat1}$ are disjoint equally long subintervals of $J_s$,
\item the left endpoints of $J_s$ and $J_{s\concat0}$ coincide,
\item the right endpoints of $J_s$ and $J_{s\concat1}$ coincide,
\item the gap between $J_{s\concat0}$ and $J_{s\concat1}$ is $2^{-\abs s}a_{\abs s}b_k$.
\end{enumerate}
Conditions~\eqref{posl} ensure that the construction is possible.
Let
$$
  C_k=\bigcap_{p\in\omega}\bigcup_{s\in 2^p}J_s
$$
be the resulting ternary set.
Let
$$
  C=\prod_{k\in\omega}C_k\subs\Tset.
$$
The first property we notice is $\lambda_\Tset(C)>1-\eps$.
Indeed, $C_k$ is omitting a set of Lebesgue measure exactly
$b_k(1+\sum_{m\in\omega}a_n)<2b_k$. Hence $\leb(C_k)>1-2b_k$. Therefore
$$
  \lambda_\Tset(C)=\prod_k(1-2b_k)\geq 1-2\sum_k b_k>1-\eps.
$$

We claim that $C$ is also nearly ultrametric. To prove it,
consider the Cantor cube $\Cset$ and provide it with the usual least difference
metric: $\rho(x,y)=2^{-\abs{x\wedge y}}$ if $x\neq y$, $\rho(x,x)=0$.
(Here and later, $x\wedge y$ denotes the initial segment common to $x$ and $y$.)
For $x\in\Cset$ denote by $\wh x$ the unique
point of $C_k$ coded by $x$.

Now consider the power $(\Cset)^\omega$ of the Cantor set $\Cset$ and equip it
with the metric
$$
  \rho_G(x,y)=\sup_{k\in\omega}r_k\rho(x_k,y_k).
$$
For $x\in(\Cset)^\omega$ let $\wh x=\seq{\wh x_k:k\in\omega}$.
The mapping $x\mapsto\wh x$ is obviously bijective. It is also Lipschitz:
Let $x,y\in(\Cset)^\omega$ and let $k\in\omega$.
Set $s=x_k\wedge y_k$. Then
$\norm{\wh x_k-\wh y_k}\leq\leb(J_s)\leq2^{-\abs s}=\rho(x_k,y_k)$.
Since this is true for every coordinate, we have
$d_G(\wh x,\wh y)\leq\rho_G(x,y)$.

We now prove that the inverse $\wh x\mapsto x$ of the map is nearly Lipschitz.
In view of Proposition~\ref{prop:nL} it is thus enough to proceed as follows:
Let $j\in\omega$ and suppose $\rho_G(x,y)=2^{-j}$. There is $n$ and $k\in I_n$ such that
$r_k\rho(x_k,y_k)=2^{-j}$. Put $s=x_k\wedge y_k$ and $p=\abs{s}$.
Since $r_k=2^{-n}$, we have $2^{-n}2^{-p}=2^{-j}$, i.e., $n+p=j$.
Estimate $d_G(\wh x,\wh y)$:
$$
  d_G(\wh x,\wh y)\geq r_k\norm{\wh x_k-\wh y_k}
  \geq r_k\dist(J_{s\concat0},J_{s\concat1})
  =r_kb_k2^{-p}a_p=2^{-j}b_ka_p.
$$
Therefore (employing definitions of $a_m$ and $b_k$)
\begin{multline*}
  \frac{\log d_G(\wh x,\wh y)}{\log\rho_G(x,y)}\leq
  \frac{j-\log a_p-\log b_k}{j} \\
  =\frac{j+2\log(p+1)+1+2\log(n+1)+2+\log G(n)-\log\eps}{j}.
\end{multline*}
Let $\tilde G(n)=\max_{i\leq n}G(i)$. It is easy to check that since $G$ is slow, so
is $\tilde G$. Since $n+p=j$, we conclude that
\begin{multline*}
  \limsup_{j\to\infty}\sup_{\rho_G(x,y)=2^{-j}}
  \frac{\log d_G(\wh x,\wh y)}{\log\rho_G(x,y)}\\
  \leq
  \limsup_{j\to\infty}
  \frac{j+4\log(j+1)+3+\log\tilde G(j)-\log\eps}{j}=1,
\end{multline*}
which is enough.
In summary, the mapping $x\mapsto \wh x$ is a nearly Lipschitz equivalence of $C$ and
$((\Cset)^\omega,\rho_G)$. Since $((\Cset)^\omega,\rho_G)$
is an ultrametric space, it follows that $C$ is nearly ultrametric.
\end{proof}

We now extend the statement to all measures on $\Tset$.
\begin{prop}\label{prop:torus}
Let $\mu$ be a finite Borel measure on $\Tset$ and $\eps>0$.
If $G$ is slow, then there is a compact nearly ultrametric set
$K\subs(\Tset,d_G)$ such that $\mu(\Tset\setminus K)<\eps$.
\end{prop}
\begin{proof}
Assume without loss of generality that $\mu(\Tset)=1$.
Let $C\subs\Tset$ be the nearly ultrametric set
from Lemma~\ref{lem:haar}.
Let $A=\{(x,y):x+y\in C\}$. By Fubini Theorem,
$(\lambda_\Tset\times\mu)(A)=\int_\Tset\lambda_\Tset(C-x)d\mu
=\lambda_\Tset(C)\mu(\Tset)$ and at the same time $(\lambda_\Tset\times\mu)(A)
=\int_\Tset\mu(C-x)d\mu$. Thus
$$
  \int_\Tset\mu(C-x)d\mu=\lambda_\Tset(C)\mu(\Tset)=\lambda_\Tset(C)>1-\eps
$$
and consequently
there is $x$ such that $\mu(C-x)>1-\eps$.
Since $C$ is nearly ultrametric and $d_G$ is a shift invariant,
it follows that $C-x$ is also nearly ultrametric.

(Observant reader may have noticed that we actually reproved the famous Christensen
characterization~\cite{MR0326293} of Haar measure zero sets in locally compact groups
in a slightly stronger setting. We could alternatively use the Christensen's theorem
as a black box.)
\end{proof}

\subsection*{Assouad's embedding revisited}

To prove Theorem~\ref{main1}, it is enough to reduce it to the case covered by
the above proposition. In order to do so we prove the following theorem.
It is akin to classical embedding theorems of Aharoni~\cite{MR0511661} and
Assouad~\cite{MR0511662,MR763553}. We will
reiterate some of the ideas of their proofs, in particular those that appear in
~\cite{MR763553}.

\begin{thm}\label{thm:embed}
For every compact non-exploding metric space $X$ there is a slow $G\in\Pset$
and a bi-Lipschitz embedding $f:X\embed(\Tset,d_G)$.
\end{thm}
\subsection*{Construction}
Let $X$ be a metric space. The symbol $B(x,r)$ denotes, as usual, a closed
ball with center $x$ and radius $r$.

Fix $\eps>0$. Suppose $S\subs X$ is a maximal
$\eps$-separated set in $X$ (i.e., $d(s,s')>\eps$ for distinct $s,s'\in S$).
Let $I$ be a finite set and $N\in\omega$ its cardinality.
Suppose that for every $x\in X$ we have
\begin{equation}\label{color}
  \abs{S\cap B(x,8\eps)}\leq N.
\end{equation}
We employ Assouad's~\cite[Lemme 2.4]{MR763553} that claims that if~\eqref{color}
holds, then there is a coloring $\chi:S\to I$ such that
if $d(s,s')\leq8\eps$, then $\chi(s)\neq\chi(s')$.
Fix such a coloring $\chi$ and define, for each $j\in I$, a function
$\phi_j:X\to\Rset$ as follows. Let $x\in X$.
The set $\chi^{-1}(j)\cap B(x,\frac32\eps)$ has at most one point -- otherwise
there would be two points in $S$ with the same color and within distance $3\eps$.
\begin{itemyze}
\item If there is a (unique) $s\in\chi^{-1}(j)\cap B(x,\frac32\eps)$, let
$\phi_j(x)=d(x,s)$,
\item otherwise let $\phi_j(x)=\frac32\eps$.
\end{itemyze}
\begin{lem}\label{phi1}
$\forall x,y\in X\ \forall j\in I\ \abs{\phi_j(x)-\phi_j(y)}\leq d(x,y)$
\end{lem}
\begin{proof}
Clearly $\phi_j(x)=\min\{\lowerd(x,\chi^{-1}(j)),\tfrac32\eps\}$, where $\lowerd$ denotes
the lower distance of a point from a set.
The inequality
$$
  \lowerd(x,\chi^{-1}(j))\leq d(x,y)+\lowerd(y,\chi^{-1}(j))
$$
is immediate from the triangle inequality. Hence
$$
  \phi_j(x)\leq\min\{d(x,y)+\lowerd(y,\chi^{-1}(j)),\tfrac32\eps\}
  \leq d(x,y)+\phi_j(y)
$$
as required.
%
%
\end{proof}

\begin{lem}\label{phi2}
If $\frac52\eps<d(x,y)\leq5\eps$, then there is $j\in I$ such that
$$\abs{\phi_j(y)-\phi_j(x)}\geq\tfrac1{10}d(x,y). $$
\end{lem}
\begin{proof}
Since $S$ is maximal $\eps$-separated, there is $s\in S\cap B(x,\eps)$.
Let $j=\chi(s)$. For every $s'\in S\cap B(y,\frac32\eps)$ we have
$$
  d(s,s')\leq d(s,x)+d(x,y)+d(y,s')\leq\eps+5\eps+\tfrac32\eps<8\eps.
$$
At the same time
$$
  d(s,s')\geq d(x,y)-d(s,x)-d(y,s')>\tfrac52\eps-\eps-\tfrac32\eps=0
$$
which proves $s\neq s'$ and thus $\chi(s')\neq\chi(s)=j$. It follows that
$\phi_j(y)=\frac32\eps$ and consequently
$$
  \abs{\phi_j(y)-\phi_j(x)}\geq\tfrac32\eps-d(x,s)\geq\tfrac12\eps=
  \tfrac1{10}\,5\eps\geq\tfrac1{10}d(x,y).
  \qedhere
$$
\end{proof}
\subsection*{Proof of Theorem~\ref{thm:embed}}
Since $X$ is compact, we may assume that its diameter is bounded by $1$.
For every $n$ let $\eps_n=2^{-n}$. Let $S_n\subs X$ be a maximal
$\eps_n$-separated set and let $G(n)\in\omega$ be a minimal number such that
$$
   \forall x\in X\quad\abs{S_n\cap B(x,8\eps_n)}\leq G(n).
$$
\begin{claim}
$G$, defined as above, is slow.
\end{claim}
Indeed, letting
$$
  Q(\eps)=\min\{Q\in\omega:
  \forall x\in X\ \exists \{x_i:i<Q\}\ B(x,\eps)\subs\bigcup_{i<Q}B(x_i,\eps/2)\}
$$
it is easy to check that $G(n)\leq Q(8\eps_n) Q(4\eps_n) Q(2\eps_n) Q(\eps_n)$
and thus
$$
  \lim_{n\to\infty}\frac{\log G(n)}{n}\leq
  \lim_{n\to\infty}
  \frac{\log Q(8\eps_n)+\log Q(4\eps_n)+\log Q(2\eps_n)+\log Q(\eps_n)}{\log\eps_n}
  =0,
$$
because $X$ is non-exploding.

Let $\omega=I_0\cup I_1\cup I_2\cup\dots$ be the partition of $\omega$
into disjoint consecutive intervals such that the length of each $I_n$ is $G(n)$,
as discussed earlier, and let $\chi_n:S_n\to I_n$ be the corresponding coloring.
For every $j\in I_n$ let $\phi_j:X\to\Rset$ be the mapping constructed above
for $\chi=\chi_n$, $S=S_n$, $\eps=\eps_n$.
It is clear that $\phi_j(X)\subs[0,\frac32\eps_n]$.
Thus the mapping $\phi$ on $X$ defined by
$\phi(x)=\seq{\frac13\phi_j(x):j\in\omega}$ maps $X$ into the cube
$\prod_{n\in\omega}[0,2^{-n-1}]^{I_n}$.
Equip this cube with the supremum metric.
Now Lemma~\ref{phi1} proves that $\phi$ is Lipschitz
and Lemma~\ref{phi2} proves that $\phi^{-1}$ is Lipschitz.
Overall, $\phi:X\embed\prod_{n\in\omega}[0,2^{-n-1}]^{I_n}$
is a bi-Lipschitz embedding. For each $n$ and $j\in I_n$ let $\psi_j$ be the
linear function that maps $[0,2^{-n-1}]$ onto $[0,\frac12]$.
The mapping $\psi=\seq{\psi_j:j\in\omega}$ is clearly an isometric embedding
of $\prod_{n\in\omega}[0,2^{-n-1}]^{I_n}$ into $(\Tset,d_G)$.
Thus $f=\psi\circ\phi$ is the required bi-Lipschitz embedding
$f:X\embed(\Tset,d_G)$.
\qed

\subsection*{Proof of Theorem~\ref{main1}}
The proof of the main theorem is now trivial: Let $\mu$ be the finite Borel measure
on $X$; we may suppose $\mu(X)=1$. Since $X$ is analytic, there is a compact
set $C\subs X$ such that $\mu(C)>1-\eps/2$. By the above embedding
theorem~\ref{thm:embed} there is a slow sequence $G\in\Pset$ and a bi-Lipschitz embedding
$C\embed(\Tset,d_G)$. Let $\nu$ be the image measure of the restriction
of $\mu$ to the set $C$. By Proposition~\ref{prop:torus} there is a
nearly ultrametric set $K\subs\Tset$ such that
$\nu(K)>\nu(f(C))-\eps/2=\mu(C)-\eps/2>1-\eps$.
The desired set is $f^{-1}(f(C)\cap K)$.
\qed

\section{Indecomposability}\label{sec:indec}

Given $n\in\omega$, if the $n$-dimensional Hausdorff measure of a metric space
$X$ is positive, then its Hausdorff dimension is at least $n$,
but not necessarily vice versa.
We will study a condition between the two, whose importance lies in the fact
that it is sufficient and also necessary for $X$
to be mapped onto an $n$-dimensional ball by a nearly Lipschitz mapping.
The main result of this section is Theorem~\ref{thm:nLsD}.
We have to go through some preliminary material first.

\subsection*{Hausdorff functions}

Recall that a right-continuous, non-decreasing function $h:[0,\infty)\to[0,\infty)$
with $h(0)=0$ is called a \emph{gauge} or a \emph{Hausdorff function}.
Recall that a gauge is \emph{doubling} if there is a constant $C$
such that $h(2r)\leq Ch(r)$ for all $r>0$.
For a gauge $h$ define
$$
  \ord h=\varliminf_{r\to0}\frac{\log h(r)}{\log r}.
$$

\subsection*{Hausdorff measures}

Recall that given a gauge $g$, the \emph{Hausdorff measure} $\hm^g$
on a metric space $X$ is defined thus: For each $\delta>0$ and $E\subs X$
set
\begin{equation}\label{hausdorff}
  \hm^g_\delta(E)=\inf\sum_ng(\diam E_n),
\end{equation}
where the infimum
is taken over all finite or countable covers
$\{E_n\}$ of $E$ by sets of diameter at most $\delta$, and put
\begin{equation*}
  \hm^g(E)=\sup_{\delta>0}\hm^g_\delta(E).
\end{equation*}
The basic properties of $\hm^g$ are well-known.
It is an outer measure and its restriction to Borel sets is
a $G_\delta$-regular Borel measure in $X$.
General references:~\cite{MR867284,MR1333890,MR0281862}.
We shall need the following theorem of Howroyd~\cite{MR1317515}
that generalizes earlier results of Besicovitch~\cite{MR0048540}
and Davies~\cite{MR0053184}.
\begin{thm}[{Howroyd~\cite{MR1317515}}]\label{BDH}
Let $X$ be an analytic metric space and $g$ a doubling gauge. If $\hm^g(X)>0$,
then there is a compact set $K\subs X$ such that $0<\hm^g(K)<\infty$.
\end{thm}
The particular case of Hausdorff measure when the gauge is given by $g(x)=x^s$
for a fixed $s>0$ is of major importance.
The corresponding Hausdorff measure is called
the \emph{$s$-dimensional Hausdorff measure} and denoted by $\hm^s$.

\subsection*{Hausdorff dimension}

The \emph{Hausdorff dimension} of a metric space $X$ is denoted and defined by
$$
  \hdim X=\sup\{s:\hm^s(X)>0\}.
$$
General references:~\cite{MR867284,MR1333890,MR0281862}.

\subsection*{Decomposability}
The following properties will turn crucial.
\begin{defn}
Let $s>0$. Say that a metric space $X$ is
\emph{$s$-decomposable} if there is
a countable cover $\{X_n\}$ of $X$ such that $\hdim X_n<s$ for all $n$,
and \emph{$s$-indecomposable} if it is not $s$-decomposable.
\end{defn}
It is easy to show that, for any metric space $X$
\begin{equation}\label{intermed1}
 \hm^s(X)>0\Rightarrow\text{$X$ is $s$-indecomposable}\Rightarrow\hdim X\geq s
\end{equation}
and also that none of the implications can be reversed. The first implication, though,
has a kind of converse in non-exploding spaces.
Note that the non-exploding property is an invariant of
nearly Lipschitz equivalence, and, as we shall see below in Lemma~\ref{dpdd}, so is
$s$-indecomposability.
\begin{thm}\label{intermed2}\label{thm:nLsD}
If $s>0$ and $X$ is a non-exploding analytic metric space, then
$X$ is $s$-indecomposable if and only if it is nearly Lipschitz equivalent to
a metric space $Z$ with $\hm^s(Z)>0$.
\end{thm}
To prove this theorem we prepare two lemmas.
\begin{lem}\label{hat}
Let $h$ be a gauge and $0<\beta\leq\ord h$. There is a gauge $\wh h$
with the following properties:
\begin{enum}
\item $\wh h$ is strictly increasing,
\item $\wh h(r)\geq h(r)+r^\beta$ on $[0,1]$,
\item $\ord\wh h=\beta$,
\item the gauge $\wh h^{1/\beta}$ is subadditive and $\ord\wh h^{1/\beta}=1$,
\item $\wh h$ is doubling, and if $\beta\leq 1$, then $\wh h$ is subadditive.
\end{enum}
\end{lem}
\begin{proof}
First set $h^*(r)=h(r)+r^\beta$. It is easy to verify that $\ord h^*=\beta$
and clearly $h^*\geq h$.
Now define
$$
\psi(r)=
  \begin{cases}
    \displaystyle\sup_{r\leq s\leq 1}s^{-\beta}h^*(s), & 0<r<1,\\
    \displaystyle h^*(1), & r\geq 1.
  \end{cases}
$$
If $\psi$ is bounded, let $\wh h(r)=\sup\psi\cdot r^\beta$.
In this case everything is trivial.
If $\psi$ is unbounded, put $\wh h(r)=r^\beta\psi(r)$ if $>0$ and $\wh h(0)=0$.
Routine calculation proves that $\wh h$ is right-continuous at $0$.

First note that obviously $h^*\leq\wh h$. We will often use it. In particular, (ii) holds.

(i) Suppose $0<r<s\leq 1$.
Routine calculation shows that since $h^*$ is right-continuous, the
supremum in the definition of $\psi$ is attained.
Therefore there is $t\geq r$ such that $\psi(r)=t^{-\beta}h^*(t)$.
If $t\geq s$, then $\psi(r)=\psi(s)$ and $\wh h(r)<\wh h(s)$ follows.
If $t<s$, then $\wh h(r)=r^\beta t^{-\beta}h^*(t)<h^*(t)\leq h^*(s)\leq\wh h(s)$.

(iii)
$\ord\wh h\leq\beta$ follows from (ii).
We show that $\ord\wh h\geq\beta$.
Suppose for the contrary that there is $\eps>0$ such that $\ord\wh h<\beta-\eps$.
Since $\ord h^*=\beta$, there is $\del>0$ such that $h^*(s)<s^{\beta-\eps}$
for all $s\leq\del$.
Since $\psi$ is unbounded, there is $s_0<\del$ such that
$s_0^{-\beta}h^*(s_0)\geq\psi(\del)$.
Since 
$\ord\wh h<\beta-\eps$,
there is $r<s_0$ such that $\psi(r)>r^{-\eps}$.
Therefore there is $s\in[r,1]$ such that $s^{-\beta}h^*(s)>r^{-\eps}$
and since $r<s_0$, it follows that $s\leq\del$. Since for every such $s$ we have
$h^*(s)<s^{\beta-\eps}$, it follows that
$r^{-\eps}<s^{-\beta}h^*(s)<s^{-\beta}s^{\beta-\eps}=s^{-\eps}$.
Hence $s<r$, a contradiction.

(iv)
Write $g=\wh h^{1/\beta}$. For $r<1$ we have
$g(r)=r\sup_{r\leq s\leq 1}(h^*)^{1/\beta}(s)/s$.
Hence $g(r)/r$ is non-increasing. Therefore
$$
  g(r+s)=r\tfrac{g(r+s)}{r+s}+s\tfrac{g(r+s)}{r+s}\leq
  r\tfrac{g(r)}{r}+s\tfrac{g(s)}{s}=g(r)+g(s).
$$
Since $\ord\wh h=\beta$, it is clear that $\ord g=1$.

(v)
Both statements can be easily derived from (iv).
\end{proof}
\begin{lem}
Let $s>0$ and $X$ be a non-exploding analytic space.
If $X$ is $s$-indecomposable, then there is a gauge $g$ with $\ord g\geq s$ and
$\hm^g(X)>0$.
\end{lem}
\begin{proof}
Since every ball in a non-exploding space $X$ is obviously totally bounded,
its closure in the completion of $X$ is compact. It follows that $X$ is contained
in a locally compact, subset of its completion.
In particular, it is a subset of a \si compact space.

We may thus employ the following \cite[Theorem 6.4]{MR0146331} of Sion and Sjerve:
if $X$ is an analytic subset of a \si compact metric space, and
$0<s_0<s_1<s_2<\dots$ is a sequence of reals, then
\begin{enumerate}
\item either there is a cover $\{X_n\}$ of $X$ such that $\hm^{s_n}(X_n)=0$
for all $n$,
\item or else there is a gauge $g$ such that $\hm^g(X)=\infty$ and
$\ord g>s_n$ for all $n$.
\end{enumerate}
So suppose $X$ is $s$-indecomposable and pick any sequence $s_n{\nearrow}s$.
Then (1) obviously fails and thus (2) yields a gauge such that $\hm^g(X)>0$
and $\ord g\geq\sup s_n=s$.
\end{proof}

\subsection*{Proof of Theorem~\ref{thm:nLsD}}
The forward implication: Suppose $X$ is $s$-indecomposable.
By the above lemma there is a gauge $g$ with $\ord g\geq s$ and $\hm^g(X)>0$.
Using Lemma~\ref{hat} we may suppose that $\ord g=s$, $g$ is strictly increasing
and $g(r)\geq r^s$.
Define a gauge $h(r)=g(r)^{1/s}$. By Lemma~\ref{hat} we may also suppose that
$h$ is strictly increasing, subadditive, $h(r)\geq r$, and $\ord h=1$.

Let $d$ be the metric of $X$. Define a new metric on $X$ by
$\rho(x,y)=h(d(x,y))$. It is indeed a metric inducing the same topology,
because $h$ is subadditive and strictly increasing.
The identity mapping $(X,\rho)\to(X,d)$ is Lipschitz, because $h(r)\geq r$.
The identity mapping $(X,d)\to(X,\rho)$ is nearly Lipschitz, because $\ord h=1$.
Thus $(X,d)$ is nearly Lipschitz equivalent to $(X,\rho)$.
Since $\rho^s(x,y)=g(d(x,y))$, we have $\hm^s(X,\rho)=\hm^g(X,d)>0$.

The reverse implication is easy:
let $Z$ be nearly Lipschitz equivalent to $X$ and $\hm^s(Z)>0$. Then
$Z$ is $s$-indecomposable by \eqref{intermed1} and,
by Lemma \ref{dpdd} below, so is $X$.
\qed

\section{Mapping non-exploding spaces onto self-similar sets and cubes}
\label{sec:maps}

We are ready to present and prove the second summit of the paper: an analytic
non-exploding space maps onto the cube $[0,1]^n$ by a nearly Lipschitz map if
and only if it is $n$-indecomposable.

We will actually prove a bit more general statement that involves self-similar sets.
The material is taken from~\cite{MR2957686}, where a \emph{self-similar set}
is defined as the attractor of an iterated function system with all functions
being contracting similarities of $\Rset^n$.
We also \emph{a priori} impose upon self-similar sets the \emph{Open Set Condition}.
All of the relevant definitions can be found in~\cite{MR2957686}.
Interested readers are referred to one of the books~\cite{MR2118797,MR867284,MR1333890}
for an overview of self-similar sets and related material.
\begin{defn}
A mapping between metric spaces $f:X\to Y$ is termed
\emph{dimension preserving} if $\hdim f(E)\leq\hdim E$ for every set
$E\subs X$. We do not \emph{a priori} impose any continuity on the mapping.
\end{defn}
It is well-known and easy to see that Lipschitz mappings are dimension preserving
(see, e.g., \cite[Lemma 6.1]{MR867284}) and it is
also easy to see that nearly Lipschitz mappings are also. It is also worth noticing
that dimension preserving maps preserve $s$-decomposability. Proof
is straightforward.
\begin{lem}\label{dpdd}
Let $f:X\to Y$ be a dimension preserving mapping onto $Y$.
If $X$ is $s$-decomposable, then so is $Y$. In particular, the conclusion holds if $f$
is nearly Lipschitz.
\end{lem}
\begin{thm}\label{mainlip}
Let $X$ be a non-exploding analytic metric space
and $S\subs\Rset^n$ a self-similar set satisfying the Open Set Condition.
Let $s=\hdim S$. The following are equivalent.
\begin{enum}
\item $X$ is $s$-indecomposable.
\item There is a nearly Lipschitz mapping $f:X\to\Rset^n$
such that $S\subs f[X]$.
\item There is a dimension preserving surjection $g:X\to S$.
\end{enum}
\end{thm}
\begin{proof}
(ii)\Implies(iii) is easy: start with the nearly Lipschitz map $f$.
Pick a single point $z\in S$ and
define $g(x)=f(x)$ if $x\in f^{-1}(S)$ and $g(x)=z$ otherwise.
Since, as pointed out above, $f$ is dimension preserving, so is $g$.

(iii)\Implies(i) is also easy: suppose that (iii) holds and yet $X$ is
$s$-decomposable. By Lemma~\ref{dpdd} $S$ is $s$-decomposable as well.
But self-similar sets with the Open Set Condition are indecomposable because
they have positive Hausdorff measure,
see, e.g., \cite{MR867284}.

The only remaining implication (i)\Implies(ii) is harder. We postpone
its proof until we gather some background material.
\renewcommand{\qed}{\relax}
\end{proof}

\begin{coro}\label{corocube}
Let $X$ be a non-exploding analytic metric space
and $n\in\omega$. The following are equivalent.
\begin{enum}
\item $X$ is $n$-indecomposable.
\item There is a nearly Lipschitz surjection $f:X\to[0,1]^n$.
\item There is a dimension preserving surjection $g:X\to[0,1]^n$.
\end{enum}
\end{coro}
\begin{proof}
This is immediate, because the cube $[0,1]^n$ is a self-similar set with
the Open Set Condition.
We only need to take care of values of $f$ that are outside $[0,1]^n$.
But since $[0,1]^n$ is convex, we may send every such point to the
closest point on the boundary of $[0,1]^n$. This
mapping is Lipschitz, which is enough, since a composition of a nearly Lipschitz map
and a Lipschitz map is nearly Lipschitz.
\end{proof}
Let us point out that the above theorem and corollary apply in particular to
Borel sets in Euclidean spaces. We present an illustration:
\begin{coro}\label{coroX}
Let $m\leq n$ be positive integers.
For every $m$-indecomposable Borel set $X\subs\Rset^n$ there is a
nearly Lipschitz surjection $f:X\to [0,1]^m$.
In particular, such a mapping exists whenever $\hm^m(X)>0$.
\end{coro}
So in particular, the Vitushkin example mentioned in the introduction
(i.e., a compact set $X\subs\Rset^2$ with positive linear measure that
cannot be mapped onto a segment by a Lipschitz map)
maps onto $[0,1]$ by a nearly Lipschitz mapping.

We are aiming towards the proof of the remaining part of Theorem~\ref{mainlip}.
We prepare a couple of notions and lemmas.

\subsection*{Monotone spaces}

At this point we make use of the notion of monotone metric space introduced
in~\cite{MR2957686}, developed in~\cite{MR2979649} and further investigated in
a number of papers, e.g., \cite{MR3158849,MR2951635,MR3439281,MR3159817,MR2822419}.
By the definition, a metric space $(X,d)$ is \emph{monotone} if
there is a linear order $<$ on $X$ and a constant $c>0$ such that
$x<y<z\Rightarrow d(x,y)<c\,d(x,z)$. The following two facts are crucial for our proof.
\begin{lem}[{\cite[Theorem 4.5, Lemma 3.2]{MR2957686}}]\label{mono1}
Let $X$ be an analytic monotone metric space and $S\subs\Rset^m$ a self-similar set.
Let $s=\hdim S$. If $\hm^s(X)>0$, then there is compact set $K\subs X$
such that $\hm^s(K)>0$ and a nearly Lipschitz mapping $g:K\to S$ onto $S$.
\end{lem}
\begin{lem}[{\cite{MR3159074,MR2979649}}]\label{mono2}
Every ultrametric space is monotone.
\end{lem}

We will also need to extend nearly Lipschitz mappings. We prove the extension lemma
in a slightly more general setting.
Let us call a mapping $f$ between metric spaces
\emph{nearly $\beta$-H\"older} if it is $\alpha$-H\"older for all $\alpha<\beta$.
This clearly extends the notion of nearly Lipschitz mapping.
Proposition~\ref{prop:nL} has a counterpart for nearly H\"older mappings:
\begin{lem}\label{holder}
A mapping $f:(X,d_X)\to(Y,d_Y)$ is nearly $\beta$-H\"older if and only if
there is a Hausdorff function $h$ such that $\ord h\geq\beta$ and
$$
  d_Y(f(x),f(y))\leq h(d_X(x,y)) \text{ for all $x,y\in X$}.
$$
\end{lem}
It is no surprise that nearly Lipschitz and nearly H\"older mappings are extendable:
\begin{lem}\label{extension}
Let $X$ be a metric space and $Y\subs X$. Let $\beta\leq1$.
Any nearly $\beta$-H\"older mapping
$f:Y\to\Rset^m$ extends to a nearly $\beta$-H\"older mapping over $X$.
\end{lem}
\begin{proof}
By Lemma~\ref{holder} and~\ref{hat} there is a subadditive Hausdorff function
$h$ such that $\ord h=\beta$ and
$\abs{f(x)-f(y)}\leq h(d(x,y))$ for all $x,y\in Y$.

It is obviously enough to prove the statement for the coordinates of $f$, so assume without
loss of generality that $f:Y\to\Rset$.
Define the extension $f^*:X\to\Rset$ by
$$
  f^*(x)=\inf_{z\in Y} f(z)+h(d(x,z)).
$$
Proving that $f^*(x)=f(x)$ for all $x\in X$ is straightforward.
We prove that $f^*$ is nearly $\beta$-H\"older. Let $x,y\in X$.
Since $h$ is subadditive, we have
$$
  f(z)+h(d(x,z))\leq f(z)+h(d(y,z))+h(d(x,y)) \text{ for all $z\in Y$.}
$$
Therefore $f^*(x)\leq f^*(y)+h(d(x,y))$ and thus
$$
  \abs{f^*(x)-f^*(y)}\leq h(d(x,y)) \text{ for all $x,y\in X$.}
$$
Apply Lemma~\ref{holder} to conclude that $f^*$ is nearly $\beta$-H\"older.
\end{proof}

\subsection*{Proof of Theorem~\ref{mainlip}(i)\Implies(ii)}

Suppose $X$ is $s$-indecomposable.
By Theorem~\ref{thm:nLsD} we may suppose that
$\hm^s(X)>0$. By the Howroyd theorem (Theorem~\ref{BDH}) there
is a finite Borel measure $\mu\leq\hm^s$ such that $\mu(X)>0$.
Apply Theorem~\ref{main1}:
there is a nearly ultrametric compact set $N\subs X$ such that $\hm^s(N)\geq\mu(N)>0$.
By Proposition~\ref{nUlt} there exists an ultrametric space $U$ and a nearly Lipschitz surjection
$\phi:N\to U$, with a Lipschitz inverse.
Since $U$ is a Lipschitz preimage of $N$, we have $\hm^s(U)>0$.

By Lemma~\ref{mono2}, $U$ is monotone, therefore Lemma~\ref{mono1} yields a compact
subset $C\subs U$ and a nearly Lipschitz mapping $g:C\to S$ onto $S$.
The composed map $g\circ\phi:N\to S$ onto $S$ is clearly nearly Lipschitz.
Now it is enough to apply Lemma~\ref{extension} to extend $g\circ\phi$ over $X$.
\qed

\section{Comments and questions}\label{sec:rem}

\subsection*{Nearly H\"older mappings}

We present a mild generalization of Theorem~\ref{mainlip}.
The dimension preserving property can be parameterized as follows:
for a mapping $f:(X,d_X)\to(Y,d_Y)$ between metric spaces define its
\emph{Hausdorff dimension} by
$$
  \hdim f=\inf\{\alpha:\hdim f(E)\leq\alpha\hdim E\text{ for every $E\subs X$}\}.
$$
It is routine to show that if $f$ is nearly $s$-H\"older, then $\hdim f\leq1/s$.

\begin{thm}\label{mainholder}
Let $X$ be a non-exploding analytic metric space and $S\subs\Rset^m$ a self-similar
set; let $s=\hdim S\geq t>0$. The following are equivalent.
\begin{enum}
\item $X$ is $t$-indecomposable.
\item There is a nearly $\frac ts$-H\"older mapping $f:X\to\Rset^m$
such that $S\subs f[X]$.
\item There is a surjection $f:X\to S$ such that $\hdim f\leq\frac st$.
\end{enum}
\end{thm}
\begin{proof}[Proof in outline]
(ii)\Implies(iii) is obvious and (iii)\Implies(i) follows easily by modification of
Lemma~\ref{dpdd}.
(i)\Implies(ii) follows the same way as in the proof of Theorem~\ref{mainlip},
one only has to insert between the spaces $U$ and $S$ an ultrametric space
$U'=(U,d^{s/t})$ where $d$ the ultrametric of $U$.
\end{proof}
\begin{coro}
Let $X$ be a non-exploding analytic metric space. Let $t\leq m$.
The following are equivalent.
\begin{enum}
\item $X$ is $t$-indecomposable.
\item There is a nearly $\frac tm$-H\"older surjection $f:X\to[0,1]^m$.
\item There is a surjection $f:X\to[0,1]^m$ such that $\hdim f\leq\frac mt$.
\end{enum}
\end{coro}

\begin{coro}[Peano curves]
For any $m>n>0$ there is a nearly $\frac nm$-H\"older Peano curve, i.e., a
surjection $p:[0,1]^n\to[0,1]^m$.
\end{coro}
This seems to have attracted some attention. E.g.,
Arnold~\cite[Problem 1988-5]{MR2078115} claims that without a proof
and asks if ``nearly'' can be dropped.
Semmes~\cite[9.1]{SEMMES} also discusses this topic.

The condition $s\geq t$ in the above theorem is not necessary,
it is only needed for the extension of the
nearly $\frac ts$-H\"older mapping over the whole space $X$, guaranteed by
Lemma~\ref{extension} for nearly $\beta$-H\"older mappings only when $\beta\leq1$.
Inspection of the proofs shows that the following remains true for any $t>0$.
Recall that we \emph{a priori} impose upon self-similar sets
the \emph{Open Set Condition}.
\begin{thm}\label{mainholder2}
Let $X$ be a non-exploding analytic metric space and $S\subs\Rset^m$ a self-similar
set; let $s=\hdim S$ and $t>0$. The following are equivalent.
\begin{enum}
\item $X$ is $t$-indecomposable.
\item There is a compact set $K\subs X$ and a nearly $\frac ts$-H\"older surjection
$f:K\to S$.
\item There is a surjection $f:X\to S$ such that $\hdim f\leq\frac st$.
\end{enum}
\end{thm}

\subsection*{Lipschitz and H\"older maps}

As mentioned in the introduction, in~\cite{MR3159074} it was proved that
any analytic metric space $X$ with $\hdim X>n$ maps onto the cube $[0,1]^n$
by a Lipschitz map. The proof builds upon ideas similar to those in the present paper,
but deviates in one detail: the constructed mapping factorizes through an interval.
Thus it is inevitably short when mappings onto disconnected
self-similar sets are under consideration.
However, there is, as indicated below, an easy remedy.
\begin{thm}\label{thm:LipMap1}
Let $X$ be an analytic metric space and $S\subs\Rset^n$ a self-similar set.
If $\hdim X>\hdim S$, then there is a compact set $C\subs X$
and a Lipschitz surjection $f:C\to S$.
\end{thm}
\begin{proof}
By the theorem of Mendel and Naor~\cite{MR3032324} quoted above
there is a compact set $C\subs X$,
a compact ultrametric space $U$ and a Lipschitz bijection $g:C\to U$ such that
$\hdim C=\hdim U>\hdim S$.
By Lemma~\ref{mono2} and~\cite[Theorem 4.7 and Lemma 3.2]{MR2957686}
there is a Lipschitz surjection $\phi:U\to S$. The required mapping is of course
$f=\phi\circ g$.
\end{proof}

%
%

\subsection*{Dimension preserving vs.~nearly Lipschitz}

As to nearly Lipschitz and dimension preserving mappings, we do not really
know anything of the relation of the notions except of the following elementary fact.
\begin{prop}\label{nldp}
A nearly Lipschitz mapping is dimension preserving and continuous.
\end{prop}
However, the theorems we proved indicate that some converse to this proposition
might hold. For example, Corollary~\ref{corocube} contains the following
information:
\emph{Let $X$ be an analytic, non-exploding space.
If there is a dimension preserving surjection $f:X\to[0,1]^m$, then
there is a nearly Lipschitz surjection $f:X\to[0,1]^m$.
}%
It is thus natural to ask:
\begin{question}
Is there any, however partial, converse to Proposition~\ref{nldp}?
\end{question}

\subsection*{Non-exploding spaces}

The notion of non-exploding space is tailored so that the proof of Theorem~\ref{main1}
works. The paramount question is, of course, whether this condition is really needed.
\begin{question}\label{q3}
Is it true that for every analytic (or, equivalently, compact) metric space $X$
and every finite Borel measure $\mu$ on $X$ there is a nearly ultrametric set
$C\subs X$ such that $\mu(C)>0$?
\end{question}
Note that if the answer were affirmative, the theorem of Balka, Darji and
Elekes~\cite{Balka2016221} quoted in the introduction (and its consequences)
would hold for all compact metric spaces.

Even a little improvement would be valuable.
Maybe there is some room for improvement.
First, there may be a more effective way of constructing a bi-Lipschitz embedding
of a compact space
into $\ell^\infty$. And maybe there is another, less regular construction of
a nearly ultrametric set in $\Tset$ with positive Haar measure that does not require
the diameters of circles to tend to zero according to a slow function $G$.
The importance of the following question is that all analytic spaces with finite
box-counting dimension satisfy the condition and thus Theorem~\ref{main1}
would hold for all analytic spaces with finite packing dimension
(since they are countable unions of sets with finite box-counting dimension).
\begin{question}\label{q4}
Let $X$ be an analytic metric space. Suppose it satisfies condition
$$
  \lim_{r\to0}\frac{\log Q(r)}{\log r}<\infty
$$
in place of condition~\eqref{def:NEXP}.
Is it true that for every finite Borel measure $\mu$ on $X$ there is
a nearly ultrametric set $C\subs X$ such that $\mu(C)>0$?
\end{question}
Let us note that an affirmative answer to either of the above questions~\ref{q3}
and~\ref{q4}
would extend the results of~\cite{MR763553}.

The other main result, Theorem~\ref{mainlip}, depends mostly on the
existence of large nearly ultrametric sets. That is another reason to study
Questions~\ref{q3} and~\ref{q4}.

\subsection*{Positive Hausdorff measure}

We know from~\cite{MR3159074} that for every analytic metric space $X$,
if $\hdim X>m$, then there is a Lipschitz surjection of $X$ onto $[0,1]^m$,
and that a non-exploding analytic $X$ is $s$-indecomposable
if and only if there is a nearly Lipschitz surjection of $X$ onto $[0,1]^m$.
The condition $\hm^m(X)>0$ is between the two.
\begin{question}
Is there a condition similar to $m$-indecomposability that characterizes
(non-exploding) compact spaces that map onto $[0,1]^m$ by a Lipschitz mapping?
\end{question}
\begin{question}
Is there a type of mapping such that the existence of such a mapping of $X$ onto $[0,1]^m$
or some similar condition
characterizes (non-exploding) compact spaces with $\hm^m(X)>0$?
\end{question}

\section*{Acknowledgments}
I would like to thank Rich\'ard Balka, M\'arton Elekes and Tam\'as Keleti
for inspiring discussions while I was visiting
Institute of Mathematics of E\"otv\"os Lor\'and University.
In particular, I am grateful to Rich\'ard Balka who pointed me to the
paper~\cite{MR0146331} and to Rich\'ard Balka and M\'arton Elekes for making me
finally write the paper. I am also grateful to the two anonymous referees
for their valuable comments.
Finally, I would like to thank the College of Wooster.
The paper was finalized during my visit thereby. In particular, my special thanks
go to Pam Pierce.

\bibliographystyle{amsplain}
\providecommand{\bysame}{\leavevmode\hbox to3em{\hrulefill}\thinspace}
\providecommand{\MR}{\relax\ifhmode\unskip\space\fi MR }
\providecommand{\MRhref}[2]{%
  \href{http://www.ams.org/mathscinet-getitem?mr=#1}{#2}
}
\providecommand{\href}[2]{#2}

\end{document}